# A NOVEL METHOD BASED ON THE TIKHONOV FUNCTIONAL FOR NON-NEGATIVE SOLUTION OF A SYSTEM OF LINEAR EQUATIONS WITH NON-NEGATIVE COEFFICIENTS

FIKS ILYA

*Institute of Applied Physics of the Russian Academy of Sciences*
*Nizhny Novgorod, 603950, 46 Ulyanov Street, Russia*
*FiksII@yandex.ru*

We propose a novel method for a solution of a system of linear equations with the non-negativity condition. The method is based on the Tikhonov functional and has better accuracy and stability than other well-known algorithms.

## 1. Introduction

A variety of applied problems can be reduced to solution of an ill-conditioned system of linear algebraic equations (SLAE). Some of them may have an additional requirement for non-negative components of an unknown vector. An unknown parameter may be, for example, concentration of the fluorophore in the diffuse fluorescence tomography (Ntziachristos *et al.* 2005, Turchin *et al.* 2008, Fiks *et al.* 2011), imaging through turbulence (Roggemann and Welsh 1996), and so on.

Thus, the problem of finding the unknown value can be formulated as follows:

$$\begin{aligned} Av &= p \\ A &\in R^{M \cdot N}, v \in R^N, p \in R^M \\ A &= \left\{ a_{ij} \geq 0, i = \overline{1,M}, j = \overline{1,N} \right\} \\ v &= \left\{ v_j \geq 0, j = \overline{1,N} \right\} \end{aligned} \quad (1)$$

Formally, this system may have no solution. In this case, we define a solution of the system as a vector which minimizes the following functional:

$$\begin{aligned} &\min F(Av, p) \\ &v_j \geq 0, j = \overline{1,N} \end{aligned} \quad (2)$$

Hereafter $v \geq 0$ means $v = \left\{ v_j \geq 0, j = \overline{1,N} \right\}$. This functional may be error $\|Av - p\|_2^2$, Tikhonov functional (Tikhonov and Arsenin 1977), and so on. The main condition is that minimization of $F(Av, p)$ must lead to minimization of the distance (in certain metrics) between $Av$ and $p$.

There exist many methods of solving SLAE. Most of them, however, do not give a solution to eq. (1), because some of them are not used for solution of non-square systems, and others do not meet the non-negativity condition. Meanwhile the non-negativity condition allows the relative condition number of eq. (1) to be reduced substantially, thus enhancing solution stability.

The first method that allowed finding minimum error under the non-negativity condition – NNLS was developed by Lawson and Hanson (Lawson and Hanson 1974). It is one of the standard methods of solving such problems but it has weak stability. Later there appeared other methods but they are less known. These include minimum entropy regularization – MER (Fletcher 1980), simultaneous multiplicative algebraic reconstruction technique – SMART (Byrne 1993), and other algorithms described in a number of works (Bertero *et al.* 1988, Hanke *et al.* 2000, Nagy and Strakos 2000, Calvetti *et al.* 2004, Jones 2006, Sun and Pitsianis 2008). Some of them have slow divergence, some weak stability, some are complicated algorithmically.

In this paper we show that with the use of the methods taking into account the non-negativity condition it is possible to obtain solutions of eq. (1) with a much higher stability than with the use of the conventional methods: matrix inversion, Tikhonov regularization, and so on. We propose a novel algorithm based on the Tikhonov functional and compare it with the methods mentioned above. The main idea underlying the proposed method is the transition from a system of linear equations to a system of quadratic equations, which is solved by a simple iteration method.

## 2. Solution stability

For analysis of solution stability we use a spectral matrix norm and the corresponding Euclidean vector norm. Then, the variation of solution $\Delta v$ caused by the variation $\Delta p$ in the right-hand side of eq. (1) is defined by

$$\frac{\|\Delta v\|_2}{\|v\|_2} \leq \frac{s_{\max}(A)}{s_{\min}(A)} \frac{\|\Delta p\|_2}{\|p\|_2}, \qquad (3)$$

where $s_{\max}(A)$ and $s_{\min}(A)$ are the corresponding maximum and minimum of singular values of matrix $A$. We can show that with the use of the condition $v \geq 0$ the upper estimate of eq. (3) may be decreased substantially, thereby a stable solution may be obtained even in the case of a high condition number of matrix $A$.

We will consider a full-rank system

$$Av = p$$
$$A \in R^{M \times N}, v \in R^N, p \in R^M. \qquad (4)$$

System (4) may have no solution in a general case, but even if there exists a solution, it may not be non-negative. We define a solution of system (4) as a solution to the problem

$$v_0 = \arg\min_{v \geq 0} \|Av - p\| \qquad (5)$$

If we define a discrepancy vector as $r_0 = Av_0 - p$, then the vector $v_0 \geq 0$ is a solution (not a unique one) of the system

$$Av = p + r_0 \qquad (6)$$

The variation $\Delta p$ in the right-side of eq. (1) will lead to variability of solution $\Delta v$ and to discrepancy vector $\Delta r$. Actually, if

$$v_\Delta = \arg\min_{v \geq 0} \|Av - (p + \Delta p)\| \tag{7}$$

then, for the obtained discrepancy vector $r_\Delta = Av_\Delta - (p + \Delta p)$, the system

$$Av = p + \Delta p + r_\Delta \tag{8}$$

has a solution $v_\Delta \geq 0$. If we denote $\Delta v = v_\Delta - v_0, \Delta r = r_\Delta - r_0$, then the difference between (8) and (6) can be written as

$$A\Delta v = \Delta p + \Delta r \tag{9}$$

We introduce the following notation:

$$m(A) \equiv \inf_{\|v\|=1, v \geq 0} \|Av\| = \inf_{\|v\| \neq 0, v \geq 0} \frac{\|Av\|}{\|v\|}$$

$$M(A) \equiv \sup_{\|v\|=1, v \geq 0} \|Av\| = \sup_{\|v\| \neq 0, v \geq 0} \frac{\|Av\|}{\|v\|}$$

The quantities $m(A)$ and $M(A)$ are analogs of singular values at $v \geq 0$.

Apparently, for the matrix $A = \{a_{ij} \geq 0, i = \overline{1,M}, j = \overline{1,N}\}$:

$$M(A) \equiv \sup_{\|x\|=1, x \geq 0} \|Ax\| = \sup_{\|x\|=1} \|Ax\| = \|A\|.$$

From the definition of $m(A)$, the following relations are valid:

$$m(A) \leq \frac{\|Av\|}{\|v\|} \tag{10}$$

$$m(A) \leq \frac{\|A(v + \Delta v)\|}{\|v + \Delta v\|}, \tag{11}$$

as $v \geq 0$ and $v + \Delta v \geq 0$.

It is obvious that, in a general case, vector $\Delta v$ does not satisfy the condition $\Delta v \geq 0$ and $\frac{\|A\Delta v\|}{\|\Delta v\|}$ may be less than $m(A)$. Using the multiplicative property of the norm for (6) we obtain

$$\|b + r_0\| \leq \|A\|\|v\| \tag{12}$$

By multiplying (11) and (12) we get

$$m(A)\|p + r_0\| \leq \frac{\|A(v + \Delta v)\|}{\|v + \Delta v\|} \|A\|\|v\| \tag{13}$$

With $\|A(v + \Delta v)\| \leq \|Av\| + \|A\Delta v\|$ and $\|v + \Delta v\| \geq \|\Delta v\| - \|v\|$ taken into account, (13) can be rewritten in the form

$$\frac{\|\Delta v\|}{\|v\|} \leq \frac{\|A\|}{m(A)} \left( \frac{\|\Delta p + \Delta r\|}{\|p + r\|} + 1 \right) + 1 \tag{14}$$

Assuming $\|\Delta p + \Delta r\| \ll \|p + r\|$ and $\|A\| \gg m(A)$, the inequality (14) takes the following form

$$\frac{\|\Delta v\|}{\|v\|} \leq \frac{\|A\|}{m(A)} + o\left(\frac{\|A\|}{m(A)}\right) \quad (15)$$

If we combine the inequalities (3) and (15), we obtain the following estimate

$$\frac{\|\Delta v\|_2}{\|v\|_2} \leq \min\left(\frac{s_{\max}(A)}{m(A)}, \frac{s_{\max}(A)}{s_{\min}(A)} \frac{\|\Delta p\|_2}{\|p\|_2}\right) \quad (16)$$

For the Euclidean norm of vectors, the magnitude of $m(A)$ can be easily calculated as shown in the lemma below.

*Lemma.* $m(A) \equiv \min_{\|x\|_2 = 1, x \geq 0} \|Ax\|_2 = \min_k \sqrt{\sum_{i=1}^{M} a_{ik}^2}$.

*Proof.*

We consider the square of the norm $\|Ax\|_2$.

$$\|Ax\|_2^2 = \sum_{i=1}^{M} (a_{i1}x_1 + a_{i2}x_2 + ... + a_{iN}x_N)^2 =$$
$$= \sum_{i=1}^{M} \left[ (a_{i1}^2 x_1^2 + a_{i2}^2 x_2^2 + ... + a_{iN}^2 x_N^2) + S_i \right],$$

where $S_i = \sum_{j_1 \neq j_2} a_{ij_1} x_{j_1} a_{ij_2} x_{j_2} \geq 0$, because $x \geq 0$ and $a_{ij} \geq 0$. Then

$$\sum_{i=1}^{M} \left[ a_{i1}^2 x_1^2 + ... + a_{iN}^2 x_N^2 + S_i \right] = x_1^2 \sum_{i=1}^{M} a_{i1}^2 + ... + x_N^2 \sum_{i=1}^{M} a_{iN}^2 + \sum_{i=1}^{M} S_i =$$
$$= \sum_{k=1}^{N} \left( x_k^2 \sum_{i=1}^{M} a_{ik}^2 \right) + \sum_{i=1}^{M} S_i. \quad (17)$$

On the assumption of $\sum_{k=1}^{N} x_k^2 = 1$ (which is a consequence of $\|x\|_2 = 1$) the expression $\sum_{k=1}^{N} \left( x_k^2 \sum_{i=1}^{M} a_{ik}^2 \right)$ has a lower limit

$$\left( x_1^2 \sum_{i=1}^{M} a_{i1}^2 + x_2^2 \sum_{i=1}^{M} a_{i2}^2 + ... + x_N^2 \sum_{i=1}^{M} a_{iN}^2 \right) \geq (x_1^2 + ... + x_N^2) \cdot \min_k \sum_{i=1}^{M} a_{ik}^2 = \min_k \sum_{i=1}^{M} a_{ik}^2 \quad (18)$$

This estimate is achieved at $x_{k^*} = 1, x_j = 0, j = \overline{1, N}, j \neq k^*$, where $k^* = \arg\min_k \sum_{i=1}^{M} a_{ik}^2$. Taking into consideration that for $x_{k^*} = 1, x_j = 0, j = \overline{1, N}, j \neq k^*$, the values $S_i = 0, i = \overline{1, M}$ we derive

$$\min_{x_1^2 + x_2^2 + ... + x_N^2 = 1} \left( x_1^2 \sum_{i=1}^{M} a_{i1}^2 + x_2^2 \sum_{i=1}^{M} a_{i2}^2 + ... + x_N^2 \sum_{i=1}^{M} a_{iN}^2 + \sum_{i=1}^{M} S_i \right) = \min_k \sum_{i=1}^{M} a_{ik}^2 \quad (19)$$

Turning to the Euclidean norm, we obtain

$$\min_{\|x\|_2 = 1, x \geq 0} \|Ax\|_2 = \min_k \sqrt{\sum_{i=1}^{M} a_{ik}^2} \quad \square$$

As a rule, $m(A)$ is significantly more than $s_{\min}(A)$ in the majority of cases.

Despite the rather crude estimate, the inequality (16) shows that accounting for the non-negativity condition greatly reduces the relative condition number of SLAE (1). Hence, from this viewpoint, it is of prime importance to use when solving SLAE algorithms giving a solution with non-negative components. Actually, this means that, if a system of equations with constraints in the form of the inequality

$$\begin{cases} Av = p \\ p \geq 0 \end{cases} \quad (20)$$

may be reduced to a system of linear equations without constraints:

$$Tv = b,$$

then, except for the degenerate cases, the conditionality of the latter system will reduce to $cond(T) \ll cond(A)$. Thereby, developing special methods for obtaining non-negative solutions is justified.

### 3. A method for obtaining non-negative solutions based on Tikhonov functional

The idea of the proposed method is to move from a constrained optimization problem to an unconstrained optimization problem and vector $v \geq 0$ can be written as $v = \{v_j = u_j^2, j = \overline{1, N}, u_j \in R\}$. The latter can be written in matrix-vector form:

$$v = D(u)u, \quad (21)$$

where $D(u) = diag(u)$ that is $D(u) = \{d_{ij} : d_{jj} = u_j, i, j = \overline{1, N}\}$. Then (1) can be rewritten as a follows:

$$AD(u)u = p, \quad (22)$$

For solution of the equation (22) we use a simple iteration

$$AD(u^{(k)})u^{(k+1)} = p \quad (23)$$

So, for $(k+1)$ the iteration matrix $AD(u^{(k)})$ is known. Based on (23) we write a Tikhonov regularization equation (Tikhonov and Arsenin 1977):

$$(D(u_\alpha^{(k)})A^T AD(u_\alpha^{(k)}) + \alpha I)u_\alpha^{(k+1)} = D(u_\alpha^{(k)})A^T p, \quad (24)$$

where $u_\alpha^{(k)}$ is a regularized solution, $I$ is an identity $N \times N$ matrix. Equation (24) can be obtained directly from (23). Matrix $D(u_\alpha^{(k)})A^T AD(u_\alpha^{(k)}) + \alpha I$ in the case $\alpha > 0$ is a nondegenerate, Hermitian and positive definite matrix and can be inverted (Tikhonov and Arsenin 1977). To solve equation (24) we will use the method of simple iterations with relaxation parameter $0 < \omega < 1$:

$$u_\alpha^{(k+1)} = \omega u_\alpha^{(k)} + (1-\omega)(D(u_\alpha^{(k)})A^T AD(u_\alpha^{(k)}) + \alpha I)^{-1} D(u_\alpha^{(k)})A^T p \quad (25)$$

The iteration procedure is stopped, if the following condition is satisfied:

$$\left\| u_\alpha^{(k+1)} - u_\alpha^{(k)} \right\| \leq \varepsilon, \quad (26)$$

where $\varepsilon$ is a specified small number.

We call this method as described the Tikhonov regularization with non-negativity constraint (TRNNC).

## 4. Testing the developed method

The developed algorithm will be studied on a system with the following matrix

$$A = \frac{1}{2}(H + \bar{H}),  \qquad (27)$$

where $H \in \mathbb{R}^{N \times N}$ is the Hilbert matrix with $N=30$ and $\bar{H}$ is the mirror matrix of matrix $H$, i.e.

$$\bar{H}_{i,j} = H_{N-i+1, N-j+1}, \; i,j = \overline{1,N} \qquad (28)$$

The matrix has the following characteristics: $s_{\max}(A) \approx 1.4, s_{\min}(A) \approx 1.6 \cdot 10^{-18}, m(A) \approx 0.2$. We will test different methods for system solution:

- General methods: matrix inversion – INV, algebraic reconstruction technique – ART (Natterer 1986, Byrne 2008), Tikhonov regularization – TR (Golub *et al.* 1997), generalized minimal residual method – GMRES (Saad and Schultz 1986).
- The methods taking into account the non-negativity condition: non-negative least squares – NNLS (Lawson and Hanson 1974), simultaneous multiplicative algebraic reconstruction technique – SMART (Byrne 2005), maximum entropy regularization – MER (Smith *et al.* 1985, Sun and Pitsianis 2008), automatically regularized non-negative AUTOREGNN (Jones, Jones 2006), modified residual norm steepest descent – MRNSD (Nagy and Strakos 2000), TRNNC.

During the test we compare the vector $v_0$ and the vector $v$ that is a solution of the equation $Av = b$, where $b = A \cdot v_0 + \varepsilon$, $\varepsilon = 10^{-10}(rnd - 0.5)\|Av_0\|_2$, with *rnd* being a random $N \times 1$ vector of uniform distribution function on the interval $(0;1)$. Accuracy of the result will be described by parameter $\rho = \frac{\|v_0 - v\|_2}{\|v_0\|_2}$. Below we present the results of 6 tests for different values of the vector $v_0$. The values of $\rho$ are indicated at the top of the corresponding figures.

*Test No.1.*
$$v_{0i} = \begin{cases} 1, i = 15 \\ 0, \text{otherwise} \end{cases}$$

*Test No.2.*
$$v_{0i} = \begin{cases} 1, 10 \leq i \leq 20 \\ 0, \text{otherwise} \end{cases}$$

*Test No.3.*
$$v_{0i} = \begin{cases} 0, 10 \leq i \leq 20 \\ 1, \text{otherwise} \end{cases}$$

| | Test No.1 | Test No.2 | Test No.3 |
|---|---|---|---|

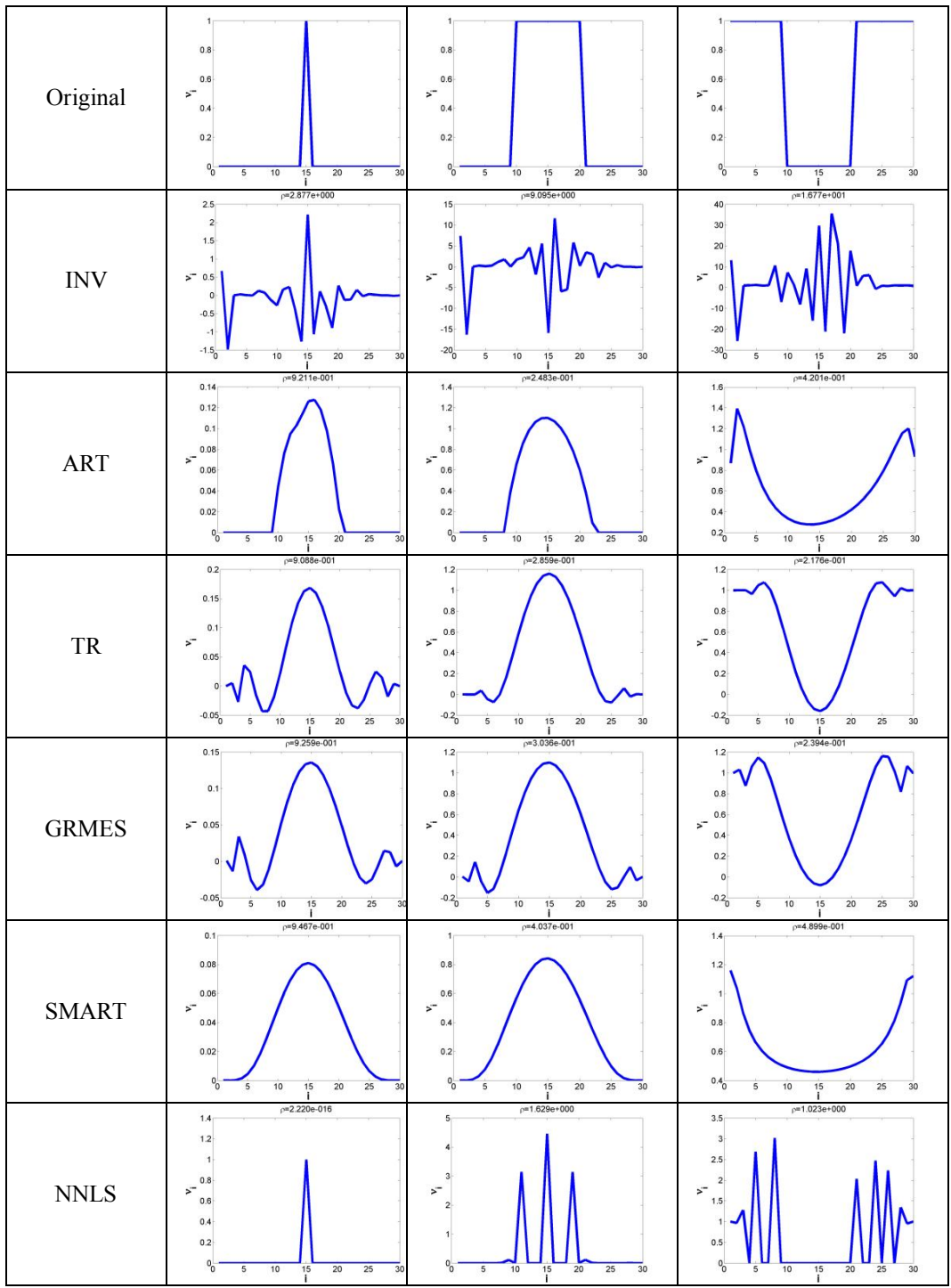

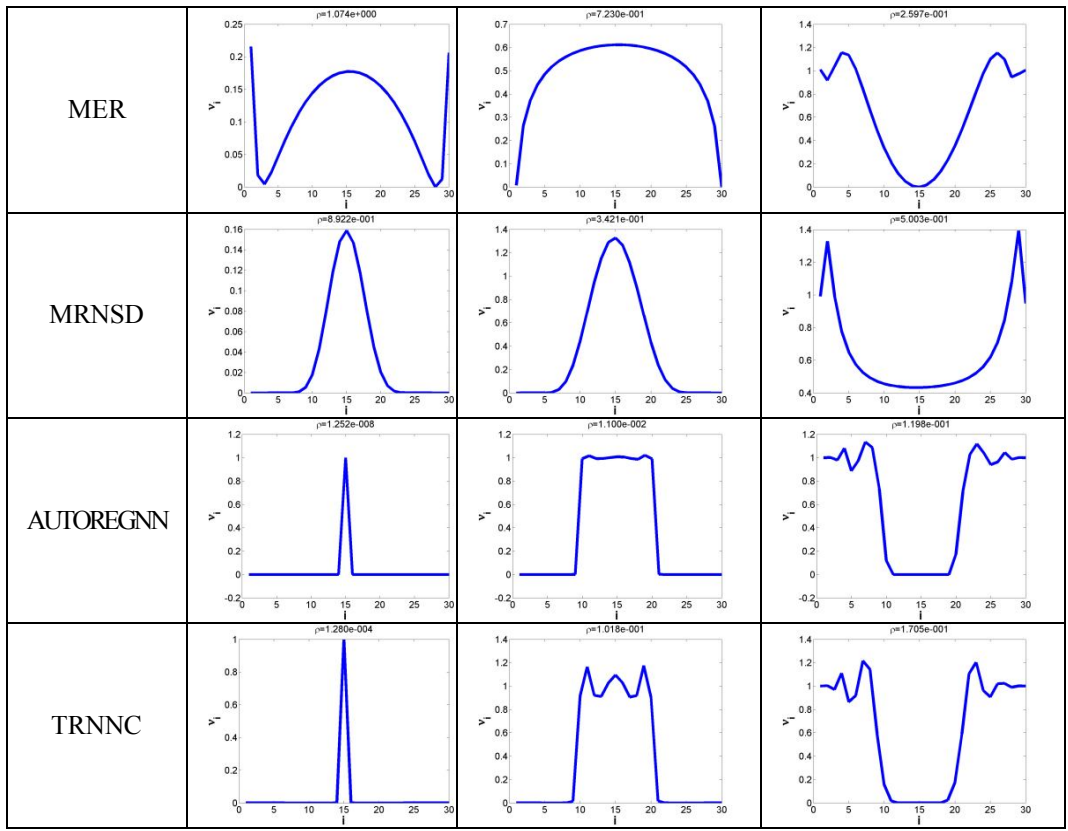

Test No.4.
$$v_{0i} = \begin{cases} \dfrac{i}{15}, 1 \le i \le 15 \\ 2 - \dfrac{i}{15}, \text{ otherwise} \end{cases}$$

Test No.5.
$$v_{0i} = e^{-\dfrac{(i-15)^2}{20}}$$

Test No.6.
$$v_{0i} = 1 + 0.5 \sin\left(\dfrac{2\pi i}{15}\right)$$

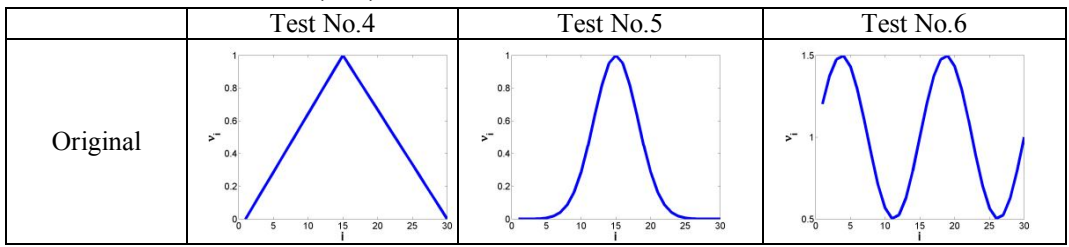

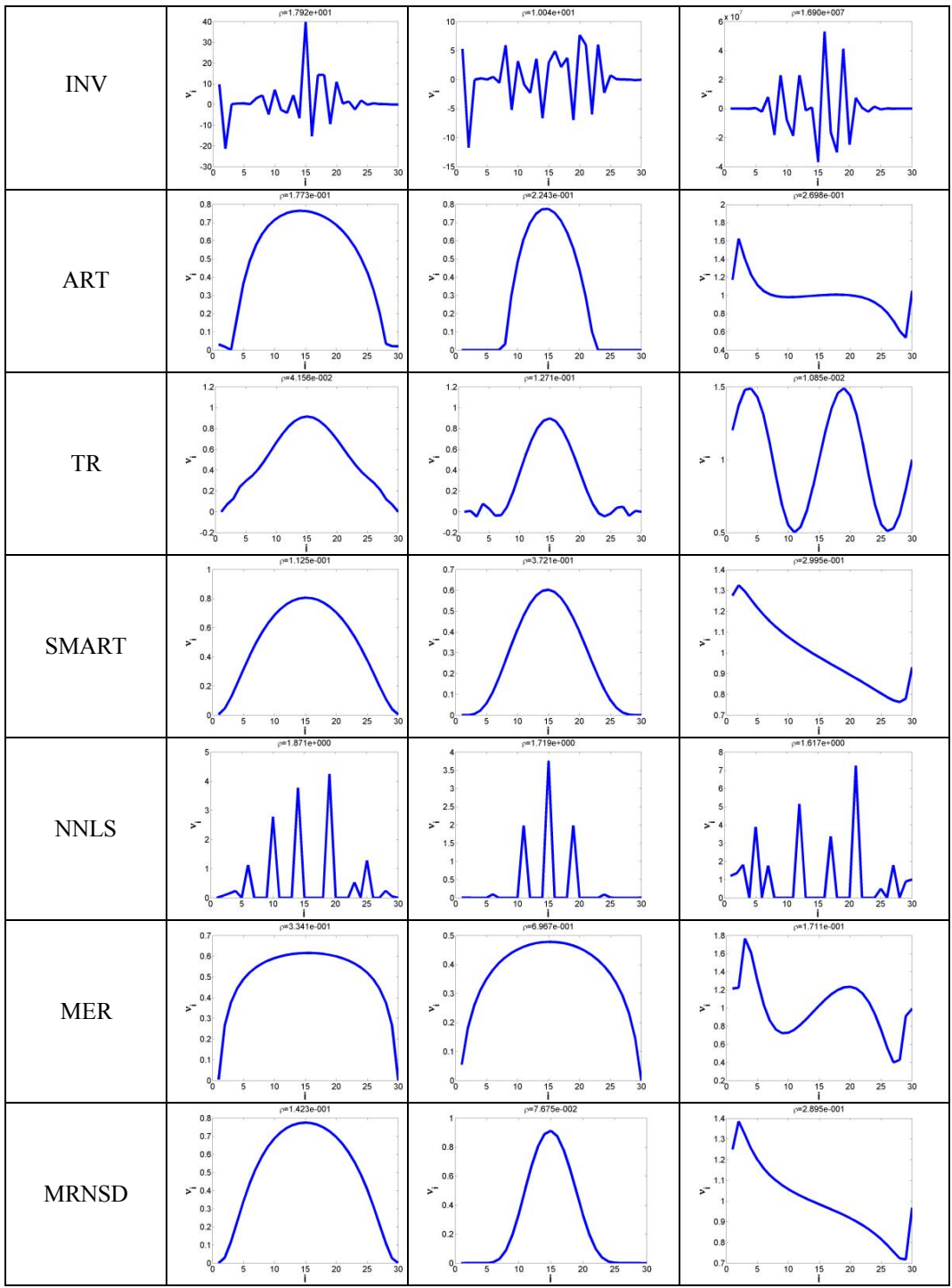

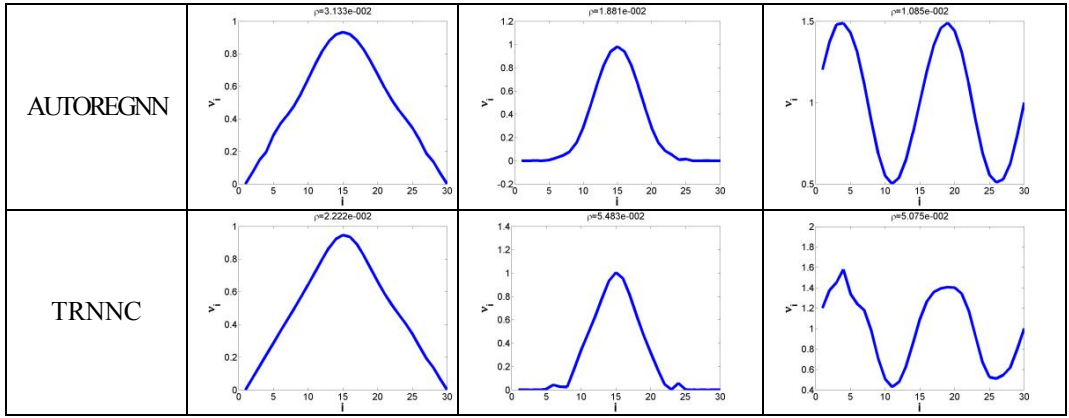

We also tested TRNNC method for systems with non-square matrixes $A_{over}$ (overdetermined matrix) and $A_{under}$ (underdetermined matrix).

$$A_{over} = A_{ij}, i = \overline{1,30}, j = \overline{1,20}$$
$$A_{under} = A_{ij}, i = \overline{1,20}, j = \overline{1,30},$$

The matrixes have the following characteristics: $s_{max}(A_{over}) \approx 1.2$, $s_{min}(A_{over}) \approx 1.5 \cdot 10^{-17}$, $m(A_{over}) \approx 0.2$, $s_{max}(A_{under}) \approx 1.2$, $s_{min}(A_{under}) \approx 1.3 \cdot 10^{-17}$, $m(A_{over}) \approx 0.15$.

Below the results of 3 tests for different values of the vector $v_0$ are presented. The initial vectors $v_0$ have distributions similar to tests No.1,2,5.

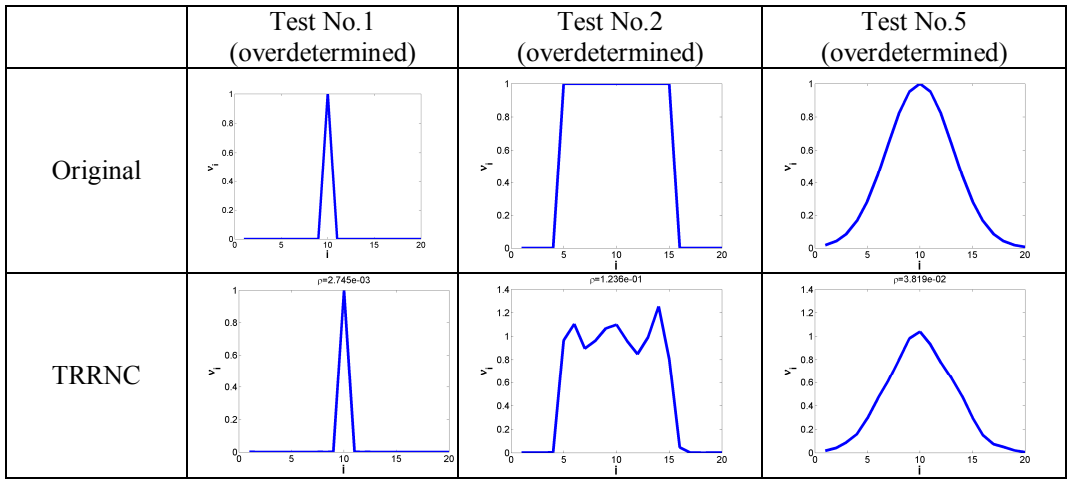

| | Test No.1 (underdetermined) | Test No.2 (underdetermined) | Test No.5 (underdetermined) |
|---|---|---|---|
| Original | 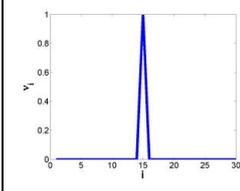 | 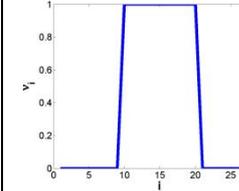 | 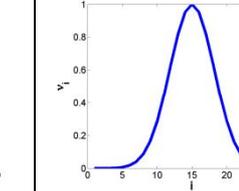 |
| TRRNC | 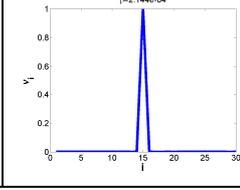 | 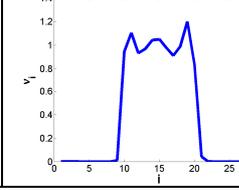 | 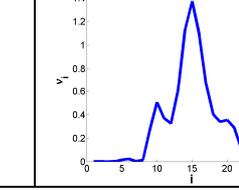 |

## 5. Results

One can see from the figures that the matrix inversion doesn't give a correct solution in any case. The best of the tested methods is AUTOREGNN. TRNNC is a quite good improvement of Tikhonov regularization. All the other methods frequently failed to correctly reconstruct the initial distribution. This is partially due to their slow rate of convergence. Note that TRNNC and AUTOREGNN are devoid of this drawback, as each iteration is, actually, finding of an exact solution of a system of equations of each iteration. A significant drawback of the AUTOREGNN method is its extremely long running time at large matrix size. It is connected with the fact that AUTOREGNN requires almost as many iterations as the size of the unknown vector (more precisely, the number of iterations is equal to the number of zero components of the solution), however, at each iteration the size of the system is reduced by 1.

The results of TRNNC method approbation for systems with non-square matrixes don't demonstrate significant difference from results with square matrixes(especially, in case of overdetermined matrix). However, it is not correct to evaluate the accuracy of underdetermined system solution, because such systems have infinite number of equivalents solutions. In this situation, it is impossible to guarantee that the method will preferably provide the initial vector opposed to other solutions.

## 6. Discussions and conclusion

The main advantage of the proposed method is search for exact solution to equation (24) at each iteration, which enhances convergence of the iteration process. On the other hand, it is a significant drawback, as the matrix $(D(u_\alpha^{(k)})A^T AD(u_\alpha^{(k)})+\alpha I)^{-1}D(u_\alpha^{(k)})A^T p$ must be calculated at each iteration, which requires additional memory and calculation time. However, for solution of the system $(D(u_\alpha^{(k)})A^T AD(u_\alpha^{(k)})+\alpha I)u_\alpha^{(k+1)} = D(u_\alpha^{(k)})A^T p$ with a symmetric positive Hermitian matrix one can use, instead of matrix inversion, iteration methods, for instance, the conjugate gradient method. This permits using less memory and, in case of a large-size $A$ matrix, reduces the computation time.

Selection of parameter $\alpha$ is not considered in the present paper. The problem is that estimating $\alpha$ is not a trivial problem. This problem is connected with changing the system to be solved at each iteration, so parameter alpha depends on $\left\|u_\alpha^{(k)}\right\|$ and should be changed at each iteration. However, calculation of this parameter at each iteration is a costly procedure. On the other hand, alpha may be rewritten, for example, as $\alpha = \beta \left\|u_\alpha^{(k)}\right\|, \beta = const$ (Cetin *et al.* 2002). In this case, the regularization parameter changes at each iteration automatically, but this formula does not provide good solution. This parameter is interpreted as a parameter which allows inverting such matrixes $D(u_\alpha^{(k)})A^T AD(u_\alpha^{(k)}) + \alpha I$, in our tests $\alpha = 10^{-13}$. More precise estimate of $\alpha$ is required considerable additional investigation.

The relaxation parameter omega needs to convergence of iteration procedure $u_\alpha^{(k+1)} = \omega u_\alpha^{(k)} + (1-\omega)T\left(u_\alpha^{(k)}\right)$. We do not consider convergence of the method either. Today, convergence of the proposed method has been proved for arbitrary initial $u_\alpha^{(0)}$ under a definite condition imposed on parameter $\alpha$. Unfortunately, this condition is a very crude estimate that is not fit for practical applications; hence it is not presented in the current study.

To conclude, a novel method for a solution of a system of linear equations with the non-negativity condition was proposed. This method is based on the Tikhonov functional and has better accuracy and stability compared to other well-known algorithms. It was shown that accounting for the non-negativity condition can greatly reduce the relative condition number of SLAE.

**Acknowledgements**

This work was financially supported by the Program "Fundamentals of Basic Studies of Nanotechnologies and Nanomaterials" of the Presidium of the Russian Academy of Sciences", the Russian Foundation for Basic Research (projects # 10–02–01109, #12-02-31361), grant of the President of the Russian Federation for support of young scientists (project # 16.120.11.1909-MK) and Measures to Attract Leading Scientists to Russian Educational Institutions program.

The author thanks Ilia Turchin for valuable suggestions and Nadezhda Krivatkina for helping with translation to English.